\input amstex
\input amsppt.sty
\magnification=\magstep1
\advance\vsize0.5cm\voffset=-1.4cm\advance\hsize1.5cm\advance\hoffset0cm
\NoBlackBoxes
\define\R{{\Bbb R}}
\define\Z{{\Bbb Z}}
\define\Q{{\Bbb Q}}

\def\pr{\mathop{\fam0 pr}}

\def\diag{\mathop{\fam0 diag}}

\def\Emb{\mathop{\fam0 Emb}}

\def\id{\mathop{\fam0 id}}
\def\dist{\mathop{\fam0 dist}}

\def\delet{\mathaccent"7017 }

\def\Cl{\mathop{\fam0 Cl}}
\def\Emb{\mathop{\fam0 Emb}}

\def\t{\widetilde}

\topmatter
\newpage
\title Embeddings of homology equivalent manifolds with boundary
\endtitle
\author D. Gon\c calves and A. Skopenkov \endauthor
\thanks
This is a new version of [GS06].
Gon\c calves is supported in part by FAPESP `Projecto  Tem\'atico Topologia Alg\'ebrica, Geom\'etrica e Differencial' number   2008/57607-6.
Skopenkov is supported in part by the Russian Foundation for Basic Research Grant No. 12-01-00748-a and
by Simons-IUM Fellowship.
We are grateful to S. Melikhov and A. Volovikov for useful discussions.
\endthanks
\address
Departamento de Matem\'atica, IME, University of S\~ao Paulo, Caixa Postal 66281,
Ag\^encia Cidade de S\~ao Paulo 05311-970, S\~ao Paulo, SP, Brasil.
e-mail: dlgoncal\@ime.usp.br \endaddress
\address  Independent University of Moscow, B. Vlasy\-ev\-skiy, 11, 119002, Moscow, Russia.
e-mail: skopenko\@mccme.ru \endaddress
\abstract
We prove a theorem on equivariant maps implying the following two corollaries:

(1) Let $N$ and $M$ be compact orientable $n$-manifolds with boundaries such
that $M\subset N$, the inclusion $M\to N$ induces an isomorphism in integral
cohomology, both $M$ and $N$ have $(n-d-1)$-dimensional spines and
$m\ge\max\{n+3,\frac{3n+2-d}2\}$.
Then the restriction-induced map $\Emb^m(N)\to\Emb^m(M)$ is bijective.
Here $\Emb^m(X)$ is the set of embeddings $X\to\R^m$ up to isotopy
(in the PL or smooth category).

(2) For a 3-manifold $N$ with boundary whose integral homology
groups are trivial and such that $N\not\cong D^3$ (or for its special
2-spine $N$)
there exists an equivariant map $\t N\to S^2$, although $N$ does not embed
into $\R^3$.

The second corollary completes the answer to the following question:
for which pairs $(m,n)$ for each $n$-polyhedron $N$ the existence of an
equivariant map $\t N\to S^{m-1}$ implies embeddability of $N$ into $\R^m$?
An answer was known for each pair $(m,n)$ except $(3,3)$ and $(3,2)$.
\endabstract
\endtopmatter

\document
This note is on the classical problem of classification of embeddings into Euclidean spaces.
For recent surveys see [Sk08, MA]; whenever possible we refer to these surveys not to original papers.
As a main tool we use the Haefliger-Wu invariant defined below.

We begin with the formulation of our main homotopy result.
Let $\t N=\{(x,y)\in N\times N\ |\ x\ne y\}$ be the {\it deleted product} of $N$.
Let $\Z_2$ act on $\t N$ and on $S^{m-1}$ by exchanging factors and antipodes, respectively.
Denote by $\pi^{m-1}_{eq}(\t N)$ be the set of equivariant maps $\t N\to S^{m-1}$ up to equivariant homotopy.
The set $\pi^{m-1}_{eq}(\t N)$ can be effectively calculated
[CF60, beginning of \S2,
Ad93, 7.1, Sk02, \S6, Sk08, \S5].
Note that $\pi^{m-1}_{eq}(\t N)=\emptyset$ for $m<n$ because $\t N\supset\t{D^n}\simeq_{eq}S^{n-1}$.

We omit $\Z$-coefficients from the notation.

\proclaim{Theorem}
Let $N$ and $M$ be compact orientable connected $n$-manifolds with non-empty boundaries such that $M\subset N$
and the inclusion $M\to N$ induces an isomorphism in cohomology.
Then the restriction-induced map $\pi^{m-1}_{eq}(\t N)\to\pi^{m-1}_{eq}(\t M)$ is bijective.
\endproclaim


This homotopy result is interesting because of the following topological corollaries.
Denote CAT = DIFF or PL.
For a CAT manifold $N$ let $\Emb_{CAT}^m(N)$ be the set of CAT embeddings $N\to\R^m$ up to CAT isotopy.
A folklore general conjecture, supported by some known results (for a survey see e.g. [RS99]) is that
$\Emb_{CAT}^m(N)$ is not changed under homology equivalence of $N$ (i.e. under a map $f:M\to N$ between manifolds inducing an isomorphism in (co)homology), in the PL case for $m\ge n+3$ and in the DIFF
case for $m\ge\frac{3n}2+2$.

\proclaim{Corollary} Let $N$ and $M$ be compact orientable $n$-manifolds
with non-empty boundaries such that $M\subset N$, the inclusion $M\to N$
induces an isomorphism in cohomology,
both $M$ and $N$ have $(n-d-1)$-dimensional spines and $m\ge\max\{n+3,\frac{3n+2-d}2\}$.
Then the restriction-induced map $\Emb_{CAT}^m(N)\to\Emb_{CAT}^m(M)$ is
bijective.
\endproclaim

Recall that

$\bullet$ a subpolyhedron $K$ of a manifold $N$ is called a {\it spine} of $N$ if $N$ is a regular neighborhood of $K$ in $N$ (or, equivalently, if $N$ collapses to $K$) [RS72].
\footnote{
We remark that for a compact connected $n$-manifold $N$ with boundary, the property of having an $(n-d-1)$-dimensional spine is close to $d$-connectedness of $(N,\partial N)$.
Indeed, for a compact connected $n$-manifold $N$ with boundary and an
$(n-d-1)$-dimensional spine, the pair $(N,\partial N)$ is homologically
$d$-connected.
On the other hand, every compact connected $n$-manifold $N$ with boundary
for which $(N,\partial N)$ is $d$-connected, $\pi_1(\partial N)=0$,
$d+3\le n$ and $(n,d)\not\in\{(5,2),(4,1)\}$, has an $(n-d-1)$-dimensional
spine [Wa64, Theorem 5.5, Ho69, Lemma 5.1 and Remark 5.2].}

$\bullet$ a closed manifold $N$ (or a pair $(N,\partial N)$) is called {\it homologically
$d$-connected}, if $N$ is connected and $H_i(N)=0$ for each $i=1,\dots,d$
(or $H_i(N,\partial N)=0$ for each $i=0,\dots,d$).

In the DIFF category the restriction $m\ge\frac{3n+2-d}2$ can be relaxed to $m\ge\frac{3n+1-d}2$.

By the Corollary, any homology ball unknots in codimension at least 3, cf. [Sc77].

The Corollary follows from the Theorem and the bijectivity of
$\alpha$-invariant [RS99, \S4, Sk02, Theorems 1.1$\alpha\partial$ and 1.3$\alpha\partial$],
which is defined as follows.
For an embedding $f:N\to\R^m$ define a map
$$\t f:\t N\to S^{m-1}\qquad\text{by}\qquad\t f(x,y)=\frac{fx-fy}{|fx-fy|}.$$
The equivariant homotopy class $\alpha(f)$ of the above-defined $\t f$ in
$\pi^{m-1}_{eq}(\t N)$ is clearly an isotopy invariant.
Thus is defined the {\it Haefliger-Wu (deleted product) invariant}
$$\alpha=\alpha_{CAT}^m(N):
\Emb\phantom{}_{CAT}^m(N)\to\pi^{m-1}_{eq}(\t N).$$

{\bf Remarks.}
(a) If in Theorem and in Corollary the inclusion-induced homomorphism $H^i(N)\to H^i(M)$
is an isomorphism only for $i\ge l>0$, then the corresponding restriction-induced maps
are bijective for $m\ge n+l$ and surjective for $m=n+l-1$.

(b) The assumption that $(N,M)$ is a codimension 0 pair is essential in the Theorem and the Corollary.
Indeed, take $N=D^p\times S^q$ and $M=S^q$.
For $m\ge3q/2+2$ we have $\#\Emb^m(S^q)=1$ while $\Emb^m(D^p\times S^q)=\pi_q(V_{m-q,p})$ can contain more than
one element (specific examples are particularly easy to find for $p=1$, when $V_{m-q,p}\simeq S^{m-q-1})$.

(c) The assumption that $N$ has boundary is not essential in the Theorem and the Corollary.
But these results are trivial for closed $N$: if $N$ is closed and $M\ne N$,
then the assumptions are never fulfilled because $H^n(N)\not \cong H^n(M)$.

(d) The conclusion of the Theorem for closed manifolds is not always fulfilled, because there are closed
manifolds non-embeddable in the same dimension as the corresponding punctured manifolds.

(e) The Theorem is clearly true for $m<n$ because both sets are empty.
We conjecture that the Theorem holds for $m=n$ and for $m=n+1$.

\smallskip
Now let us present motivation for the second corollary of the Theorem.
From the construction of the map $\t f$ above it follows that

(*) if $N$ embeds into $\R^m$, then there exists an equivariant map $\t N\to S^{m-1}$.

The existence of an equivariant map $\t N\to S^{m-1}$ can be checked for
many cases [CF60, beginning of \S2,
Ad93, 7.1, Sk08, \S5].
Thus if a converse to (*) is true, the embedding problem is reduced to
a manageable (although not trivial) algebraic problem.
So in 1960s there appeared a problem to find conditions under which the converse to (*) is true.
The converse for (*) was known to be

{\it true for an $n$-polyhedron $N$ and $2m\ge3n+3$ or $m=2n=2$} [RS99, \S4, Sk08, \S5], cf. [Sk98, Theorem 1.3];

{\it false for each pair $(m,n)$ such that $\max\{4,n\}\le m\le\frac{3n}2+1$
and some $n$-polyhedron $N$} [RS99, \S4, Sk08, \S5], cf. [Sk98, Example 1.4].

In the only remaining cases $m=3$ and $n\in\{2,3\}$ it was unknown if the converse to (*) is true.
The counterexamples to the converse of (*) for $m=n\ge4$ and $m=n+1\ge4$
[MS67, Hu88] cannot be directly extended to $m=3$ because they used $m$-dimensional
contractible manifold distinct from the $m$-ball, which apparently does not exist for $m=3$.

Recall that a {\it homology $n$-ball} is an $n$-manifold with boundary whose homology groups are the same as
those of the $n$-ball.
A {\it special spine} is defined e.g. in [Ca65].

\proclaim{Proposition} The converse to (*) is false in the cases $m=3$ and $n\in\{2,3\}$:
if $N$ is either a non-trivial homology ball or a special spine of a non-trivial homology ball,
then $N$ does not embed into $\R^3$ but there exists an equivariant map $\t N\to S^2$.
\endproclaim

\demo{Proof} The non-embeddability follows because if a special spine of a
homology ball $N$ embeds into $\R^3$, then the regular neighborhood in $\R^3$
of this spine is homeomorphic to $N$ [Ca65], which contradicts to the
non-triviality of $N$.

It suffices to prove the existence of an equivariant map $\t N\to S^2$ for a homology 3-ball $N$.
\footnote{This existence follows from the (unproved) case $m=n=3$ of the Theorem because an inclusion
of the standard ball into $N$ induces isomorphisms in cohomology.}
Analogously to [Ad93, end of \S7.1] (or by Lemma 2 below) it suffices to prove that $H^i(\t N)=0$ for each $i\ge3$.
We prove this for $i=3$; the proof for each $i\ge4$ is analogous.
Let $\Delta$ be the interior of a closed regular neighborhood in $N\times N$ of the diagonal.
Then
$$H^3(\t N)
\cong H^3(N\times N-\Delta)
\cong H_3(N\times N-\Delta,\partial(N\times N-\Delta))\cong$$
$$\cong H_3(N\times N,\Cl\Delta\cup\partial(N\times N))
\cong H_2(\Cl\Delta\cup\partial(N\times N))=0,\quad\text{where}$$
\quad
$\bullet$ the first isomorphism follows because $N\times N- \Delta$ is a
deformation retract of $\t N$,

$\bullet$ the second one by Lefschetz duality (recall
that $\t N$ is orientable if $N$ is a homology ball),

$\bullet$ the third one by excision,

$\bullet$ and the fourth one by exact sequence of pair.

Using the Mayer-Vietoris sequence for
$$\partial(N\times N)=N\times\partial N
\bigcup\limits_{\partial N\times\partial N} \partial N\times N$$
and noting that $\partial N\cong S^2$, we prove that
$H_2(\partial(N\times N))=0$.
Using the Mayer-Vietoris sequence for $\Cl\Delta\cup\partial(N\times N)$ and noting
that $\Delta\simeq N$ and $\Cl\Delta\cap\partial(N\times N)$ is a regular neighborhood in
$N\times N$ of the diagonal of $\partial N$, i.e. is homotopy equivalent to
$\partial N\cong S^2$, we prove the last isomorphism.
\qed
\footnote{Another proof of the Proposition could possibly be obtained by using the fact
that for the homology 3-ball $N$, which is a punctured boundary of the Mazur
4-manifold, there exists an equivariant map $\Sigma\t N\to S^3$ [MRS03].
The obstruction to equivariant desuspension of this map on $\t P$ (where $P$ is the special spine of $N$) apparently lies in $H^4(\t P)$, which group is trivial because $P$ is acyclic [We68].}
\enddemo

Manifolds $\t N$ and $\t M$ are homotopy equivalent to $(2n-1)$-dimensional CW complexes.
Hence the Theorem follows by the cases $l=0$ of Lemmas 1 and 2, see below.
Remark (a) follows by (the general cases of) Lemma 1 and 2.

\smallskip
{\bf Lemma 1.} {\it Let $N$ and $M$ be compact orientable connected
$n$-manifolds with non-empty boundaries such that $M\subset N$ and the
inclusion induces an isomorphism $H^i(N)\to H^i(M)$ for $i\ge l$.
Then $H^i(\t N,\t M)=0$ for each $i\ge n+l$. }

\smallskip
{\bf Lemma 2.} [BG71, 3.2]
{\it Suppose $X$, $Y$ are finite connected CW-complexes with free involutions, $f:X\to Y$ is an equivariant map and
$l$ is a non-negative integer.
If  $f^*:H^i(Y)\to H^i(X)$ is an isomorphism for each $i>l$ and is onto for $i=l$, then

$(d_l)$
$f^{\sharp}:\pi^i_{eq}(Y)\to \pi^i_{eq}(X)$ is a 1-1 correspondence for $i>l$ and is onto for $i=l$.}

\smallskip
We give a proof of Lemma 2 (which was not presented in [BG71]) using standard argument and
following [HH62, pp. 236-237], cf. [Me09, Proof of Lemma 8.1].
Lemma 2 was used in the previous version [GS06] of this paper;
the proof was essentially presented there but contains mistakes which are corrected here.

\smallskip
{\it Proof of Lemma 1 for $l=0$.}
Let $N_0$ and $M_0$ be the interiors of $N$ and $M$, respectively.
It suffices to prove Lemma 1 for $N$ and $M$ replaced by $N_0$ and $M_0$.

(Indeed, the collaring theorem for the boundary of a manifold states
that {\it there is a neighborhood of  $\partial M$ in $M$ which is homeomorphic
to the product $\partial M \times [0,1)$ so that $\partial M \times \{0\}$
is mapped homeomorphically  to the boundary.}
Therefore there is an embedding $\phi:N \to N_0$ which is a homotopy
inverse of the inclusion $N_0 \to N$.
Analogously $\phi \times \phi: \t N \to \t N_0$ is a homotopy inverse of the inclusion $\t N_0 \to \t N$.
Same observations hold for $N$ replaced by $M$.
So it suffices to prove Lemma 1 for $N$ and $M$ replaced by $N_0$ and $M_0$.)

Let $x_0\in M_0 \subset N_0$ be a base point for $M_0$ and $N_0$.
Consider the following mapping of bundles (which are given by projections
onto the first factor):
$$\minCDarrowwidth{12pt}\CD M_0-x_0 @>>  > \t M_0 @>>  > M_0\\
@VV \subset V @VV \subset V @VV \subset V \\
N_0- x_0 @>>  > \t N_0 @>> > N_0 \endCD$$
The action of $\pi_1(M_0)$ in  the cohomology $H^i(M_0-x_0)$ of the fiber is trivial for each $i$.

(Indeed,  this follows for $i=n$ because $H^n(M_0-x_0)=0$ and for $i<n-1$ because $H^i(M_0-x_0)\cong H^i(M_0)$ and
the bundle is the restriction of the trivial bundle $M_0 \times M_0 \to M_0$.
For $i=n-1$ we have $M_0-x_0\simeq M_0\vee S^{n-1}$, so $ H^{n-1}(M_0-x_0) \cong  H^{n-1}(M_0)\oplus\Z$.
The action of an element $\alpha \in \pi_1(M_0)$ is given by the identity
on the first summand and multiplication by the sign of the loop on $\Z$.
Since $M_0$ is orientable, the action is identical.)

The same holds for the second bundle, where $M$ is replaced by $N$.

By excision the inclusion of the pairs $(M_0,M_0-x_0)\to(N_0,N_0-x_0)$ induces an isomorphism in cohomology.

\smallskip
{\it Proof of Lemma 1: completion for $l=0$.}
Applying 5-lemma for the inclusion-induced mapping of exact sequences of these pairs we obtain that the inclusion $M_0-x_0 \to N_0-x_0$ induces an isomorphism in cohomology.
Hence using the triviality of the action and the Universal Coefficients Theorem we obtain that
the restriction induces an isomorphism
$$r:H^p(N_0;H^q(N_0-x_0))\to H^p(M_0;H^q(M_0-x_0))\quad\text{for each}\quad p,q.$$
This $r$ is a homomorphism of the $E_2$-terms of the Leray-Serre cohomology spectral sequences of the above bundles.
By the Zeeman Comparison Theorem of spectral sequences [Ze57],
the restriction $H^i(\t N_0)\to H^i(\t M_0)$ is an isomorphism for each $i$.
This implies Lemma 1.  \qed
\footnote{A statement on cohomology of compact manifolds should have a proof involving
only cohomology of compact manifolds (recall that we may assume that
$\t N=\t N_\varepsilon$ is compact).
The above proof has such an interpretation in terms of only compact spaces.
Lemma 1 can also be proved analogously to proof of the Proposition above.}

\smallskip
{\it Proof of Lemma 1: completion for the general case.}
Applying the 5-lemma for the inclusion-induced mapping of exact sequences of these pairs we obtain that the inclusion $M_0-x_0 \to N_0-x_0$ induces an isomorphism in $H^i$ for $i\geq l$.
Hence using the triviality of the action and the Universal Coefficients Theorem we obtain that
the restriction induces an isomorphism
$$r:H^p(N_0;H^q(N_0-x_0))\to H^p(M_0;H^q(M_0-x_0))\quad\text{for}\quad p+q\ge n+l-1.$$
Hence $r$ is an isomorphism of for $p+q\ge n+l$ and an epimorphism for $p+q=n+l-1$.
This $r$ is a homomorphism of the $E_2$-terms of the Leray-Serre cohomology spectral sequences of the above bundles.
Now using standard argument of homological algebra as in the Zeeman Comparison Theorem of spectral sequences [Ze57]
we obtain that the restriction-induced homomorphism between  $E^{p,q}_r$ terms is an isomorphism for $p+q\ge n+l$
and an epimorphism for $p+q=n+l-1$.
Since $E_{n-l}=E_{n-l+1}=...=E_\infty$, the restriction induces on $E_\infty$ terms an isomorphism
for $p+q\ge n+l$ and an epimorphism for $p+q=n+l-1$.
Hence the restriction $H^i(\t N_0)\to H^i(\t M_0)$ is an isomorphism for each $i\ge n+l$ and an epimorphism for $i=n+l-1$.
Therefore by the exact sequence of pair $H^i(\t N_0,\t M_0)=0$ for each $i\ge n+l$.
Hence $H^i(\t N,\t M)=0$ for each $i\ge n+l$.
\qed


\smallskip
{\it Proof of Lemma 2.}
We may assume that $f:X\to Y$ is an inclusion.
Consider the following assertion:

$(c_l)$ {\it $H^i(Y',X';G_\varphi)=0$ for each $i>l$, finitely-generated abelian group $G$, involution
$\varphi:G\to G$ and local coefficient system $G_\varphi$ associated to $\varphi$ and double cover $(Y,X)\to(Y',X')$.}

(Local coefficient system $G_\varphi$ is defined by the following action of $\pi_1(Y')$ on $G$.
Take a representative $\alpha':[0,1]\to Y'$, $\alpha'(0)=\alpha'(1)$, of $[\alpha']\in\pi_1(Y')$.
Take a lift $\alpha:[0,1]\to Y$ of $\alpha'$.
If $\alpha(0)=\alpha(1)$, then $[\alpha']$ acts identically on $G$.
If $\alpha(0)\ne\alpha(1)$, then $[\alpha']$ acts as $\varphi$.
Clearly, this action is well-defined.)

Since $Y$ is finite-dimensional,
\footnote{It would be interesting to know if Lemma 2 holds for infinite-dimensional complexes.
Note that it does hold for infinite-dimensional complexes $S^{l-1}\to S^\infty$.}
$(c_l)$ holds for large enough $l$.
Consider the following part of the Smith-Richardson-Thom-Gysin sequences
associated to the double cover $(Y,X)\to(Y',X')$
 (see the Smith-Richardson-Thom-Gysin Sequence Theorem  below):
$$0=H^i(Y,X;G)\to H^i(Y',X';G_\varphi)\to H^{i+1}(Y',X';G_{-\varphi}).$$
By the hypothesis of Lemma 2 $H^i(Y,X)=0$ for each $i>l$.
So by the Universal Coefficients Formula $H^i(Y,X;G)=0$ for each $i>l$.
Then by downward induction on $l$ we obtain $(c_l)$.

Denote by $a$ the involution on $\pi_k(S^i)$ induced by the antipodal involution on $S^i$.
\footnote{Note that $a=\id$ for $i$ odd and $a=-\id$ for $i$ even and $k\le 2i-2$.}
The obstructions to extension to $Y$ of an equivariant map $X\to S^i$, and to homotopy uniqueness of such an extension, assume values in $H^{k+1}(Y',X';\pi_k(S^i)_a)$ and $H^k(Y',X';\pi_k(S^i)_a)$.
\footnote{This can be deduced either from obstruction theory for extension of maps with non-simply-connected range
$\R P^\infty$ [HW60]
or analogously to [CF60, beginning of \S2, Ad93, 7.1] as follows.
Denote by $t$ the involution on $Y$ and its restriction to $X$.
Define a bundle
$g:\dfrac{Y\times S^i}{(x,s)\sim(tx,-s)}\overset{S^i}\to\to Y'\quad\text{by}\quad g[(x,s)]=[x].$
Equivariant maps $Y\to S^i$ up to equivariant homotopy are in 1--1 correspondence with cross-sections
of $g$ up to equivalence.
So the required obstructions are obstructions to
\newline
(*) extendability of a section on $X'$ to a section on $Y'$ for each $i\ge l$, and to
\newline
(**) uniqueness of such an extension (up to equivalence) for $i>l$.
\newline
The action of $\pi_1(Y')$ on homotopy groups of the fiber $S^i$ gives rise to local coefficient system $\pi_k(S^i)_a$.}
These groups are trivial for $k<i$ because $\pi_k(S^i)=0$, and for $k\ge i>l$ by $(c_l)$.
So $(d_l)$ holds.
\qed

\smallskip
For a reader's convenience we present the following slight and possibly known
extension of the Smith-Richardson-Thom-Gysin sequence.
Cf. [Me09, arxiv v4, Remark 2.3 and p.9, lines 14-25].

\proclaim{Smith-Richardson-Thom-Gysin Sequence Theorem}
Let $X'$ be a connected space, $X\to X'$ a double covering and $G$ a module with an involution $\varphi$.
Consider the local coefficient system $G_\varphi$ on $X'$ associated to the double covering and $\varphi$.
Then there is a long exact sequence
$$\dots\to H^{p-1}(X';G_\varphi) \to H^p(X';G_{-\varphi}) \to H^p(X;G)\to H^p(X';G_\varphi) \to H^{p+1}(X';G_{-\varphi})\to\dots$$
If 2 is invertible in $G$ (in particular, if either $G=\Q$ or $G=\Z_p$ for $p$ an odd prime),
then we have splittable short exact sequence
$$0 \to H^p(X';G_{-\varphi}) \to H^p(X;G) \to H^p(X';G_{\varphi})\to0\quad\text{so that}$$
$$H^p(X;G) \cong H^p(X';G_{-\varphi}) \oplus H^p(X';G_\varphi).$$
If $G=\Z$ and $\varphi=\id$, then we get long exact sequence
$$\dots \to H^{p-1}(X') \to H^p(X';\Z_{-\id}) \to H^p(X)\to H^p(X') \to H^{p+1}(X';\Z_{-\id}) \to \dots$$
If $G=\Z$ and $\varphi=-\id$, then we get long exact sequence
$$\dots \to H^{p-1}(X';\Z_{-\id}) \to H^p(X') \to H^p(X)\to H^p(X';\Z_{-\id}) \to H^{p+1}(X') \to \dots$$
\endproclaim


{\it Proof.}
Consider the fibration $F\to X\to X'$ which is the double covering, where $F$ is a two-points set.
For the spectral sequence with local coefficients [Si97, Theorem 2.9] we have $E_2^{p,q}=H^p(X',H^q(F;G)_\tau)$,
where the coefficients are twisted according to double cover $X\to X'$ and the following involution $\tau$ of $H^q(F;G)$:

$\bullet$ $H^q(F;G)=0$ and $\tau$ is trivial for $q>0$, and

$\bullet$ $H^0(F;G)\cong G\oplus G$ and $\tau(a,b):=(\varphi(b),\varphi(a))$.

Then the spectral sequence contains at most one non-vanishing line.
Hence
\footnote{The twisting of $H^q(F;G)$ is as required by [Si97, 2.7].
Note that [Si97, 2.8] is not required for the statement of [Si97, Theorem 2.9] (but is required for the proof).
Note that the purpose of [Si97, Theorem 2.9] was to calculate cohomology of the total space of a fibration
not with any non-twisted coefficient system but with the twisted coefficient system coming from
a twisted coefficient system in the base.
\newline
The isomorphism $H^p(X;G)\cong H^p(X',(G\oplus G)_\tau)$
has two simpler proofs not involving spectral sequences.
According to S. Melikhov, it follows easily from definitions, as explained in [Ha, Example 3.H.2]
(the case of arbitrary $\varphi$ follows from the case $\varphi=\id$ because the involution $(a,b)\mapsto(\varphi(b),\varphi(a))$ is obtained from the involution $(a,b)\mapsto(b,a)$ by an automorphism of $G\oplus G$ [Br82, Corollary III.5.7]),
or, alternatively, is a special case of the Vietoris Mapping Theorem, [Br97,  Theorem 11.1].
}
$$H^p(X;G) \cong E_\infty^{p,0} \cong E_2^{p,0}\cong H^p(X',(G\oplus G)_\tau).$$
Let $H=\{(m,-m)\in G\oplus G\ |\ m\in G\}$.
 We have $\tau(m,-m)=(\varphi(-m),\varphi(m))=(-\varphi(m),\varphi(m))$.
Hence $\tau(H)=H$ and $(H,\tau|_H)\cong (G,-\varphi)$.
Then $(G\oplus G)/H$ has `the quotient' involution $\tau/H$.
Clearly, $((G\oplus G)/H,\tau/H)\cong (G,\varphi)$.
Now the first part of the theorem follows from the cohomological long exact sequence associated with the short
exact sequence of twisted coefficients $(H,\tau|_H)\to (G\oplus G,\tau)\to ((G\oplus G)/H,\tau/H)$.

The `further' part where 2 is invertible follows from the fact that the above short exact sequence splits:
the homomorphism $s:(G\oplus G)/H \cong G \to G\oplus G$ defined by $s(m)=(m/2,m/2)$ respects involutions and
is a splitting.
The `further' part where $G=\Z$ is clear.
\qed



\enddocument